\documentclass[11pt]{article}
 \usepackage{amsmath,amssymb,amscd}
\newcommand{\pf}[1]{\trivlist \item[\hskip\labelsep\it #1\ ]}
\newcommand{\varpf}[1]{\trivlist \item[\hskip\labelsep\sc #1:]}
\newcommand{\qedbox}{$\rlap{$\sqcap$}\sqcup$}
\newcommand{\qed}{\qquad \qedbox \endtrivlist}
\newcommand{\varqed}{\hfill \rule{0.6em}{0.6em} \endtrivlist}
\newenvironment{proof}{\pf{Proof}}{\qed}

\newenvironment{items}{
  \begin{enumerate} 
                    
  }{\end{enumerate}}

\newenvironment{keywords}{\noindent\small {\it Keywords\/}:}{\vskip 4pt}
\newenvironment{classification}{\noindent\small 2000 {\it Mathematics Subject
Classification\/}:}{\vskip 12pt}

\newcommand{\cstar}{{C^\ast}}

\newcommand{\id}{{\mathrm{id}}}

\newcommand{\A}{{\mathfrak A}}

\newcommand{\VN}{\operatorname{VN}}

\newtheorem{theorem}{Theorem}

\newtheorem{corollary}{Corollary}

\title{The operator amenability of uniform algebras}
\author{\it Volker Runde}
\date{}
\begin{document}
\maketitle
\begin{abstract}
We prove a quantized version of a theorem by M.\ V.\ She\u{\i}nberg: A uniform algebra equipped with its canonical, i.e.\ minimal, operator space structure is operator amenable if and only if it is a 
commutative $\cstar$-algebra. 
\end{abstract}
\begin{keywords}
uniform algebras, amenable Banach algebras, operator amenability, minimal operator space.
\end{keywords}
\begin{classification}
46H20, 46H25, 46J10 (primary), 46J40, 47L25.
\end{classification}
In \cite{Joh1}, B.\ E.\ Johnson introduced the notion of an amenable Banach algebra. It is an active area of research to determine, for a particular class of Banach algebras, which algebras in that
class are the amenable ones. For example, Johnson himself proved that a locally compact group $G$ is amenable if and only if its group algebra $L^1(G)$ is amenable (this characterization motivates the
choice of terminology). The characterization of the amenable $\cstar$-algebras is a deep result due to several authors (see \cite[Chapter 6]{Run} for a self-contained exposition): A $\cstar$-algebra is
amenable precisely when it is nuclear. The amenability of algebras of compact operators on a Banach space $E$ is related to certain approximation properties of $E$ (\cite{GJW}). In \cite{She},
M.\ V.\ She\u{\i}nberg showed that a uniform Banach algebra is amenable if and only if it is already a commutative $\cstar$-algebra. In this note, we prove a quantized version of She\u{\i}nberg's
theorem (and thus answer \cite[Problem 31]{Run})
\par
Our reference for the theory of operator spaces is \cite{ER}, whose notation we adopt.
\par
A Banach algebra which is also an operator space is called a completely contractive Banach algebra if multiplication is a completely contractive bilinear map. For any Banach algebra $\A$,
the maximal operator space $\max \A$ is a completely contractive Banach algebra. In \cite{Rua}, Z.-J.\ Ruan introduced a variant of Johnson's definition of amenability for completely contractive
Banach algebras called {\it operator amenability\/} (see \cite[Section 16.1]{ER} and \cite[Chapter 7]{Run}). A Banach algebra $\A$ is amenable in the sense of \cite{Joh1} if and only if
$\max \A$ is operator amenable (\cite[Proposition 16.1.5]{ER}). Nevertheless, operator amenability is generally a much weaker condition than amenability. For any locally compact group $G$,
the Forier algebra $A(G)$ carries a natural operator space structure as the predual of $\VN(G)$. In \cite{Rua}, Ruan showed that $A(G)$ --- equipped with this natural operator space structure ---
is operator amenable if and only if $G$ is amenable; on the other hand, there are even compact grups $G$ for which $A(G)$ fails to be amenable (\cite{Joh2}).
\par
Let $\A$ be a uniform algebra, i.e.\ a closed subalgebra of a commutative $\cstar$-algebra. The canonical operator space structure $\A$ inherits from this $\cstar$-algebra turns it into a
completely contractive Banach algebra. By \cite[Proposition 3.3.1]{ER}, this canonical operator space structure is just $\min \A$.
\par
We have the following operator analogue of She\u{\i}nberg's theorem:
\begin{theorem}
Let $\A$ be a uniform algebra such that $\min \A$ is operator amenable. Then $\A$ is a commutative $\cstar$-algebra.
\end{theorem}
\begin{proof}
Without loss of generality suppose that $\A$ is unital with compact character space $\Omega$. We assume towards a contradiction that $\A \subsetneq {\cal C}(\Omega)$. Combining
the Hahn--Banach theorem with the Riesz representation theorem, we obtain a complex Borel measure $\mu \neq 0$ on $\Omega$ such that
\[
  \int_\Omega f \, d\mu = 0 \qquad (f \in \A).
\]
Let $H := L^2(|\mu|)$. The canonical representation of ${\cal C}(\Omega)$ on $H$ as multiplication operators turns $H$ into a left Banach ${\cal C}(\Omega)$-module. Let $H_c$
denote $H$ with its column space structure (see \cite[p.\ 54]{ER}). Then $H_c$ is a left operator ${\cal C}(\Omega)$-module.
\par
Let $K$ denote the closure of $\A$ in $H$; clearly, $K$, is an $\A$-submodule of $H$. Trivially, $K$ is complemented in $H$, and by \cite[Theorem 3.4.1]{ER}, $K_c$ is completely complemented in $H_c$, 
i.e.\ the short exact sequence 
\begin{equation} \tag{\mbox{$\ast$}}
  \{ 0 \} \to K_c \to H_c \to H_c/ K_c \to \{ 0 \}
\end{equation}
of left operator $\A$-bimodules is admissible. Since $\min \A$ is operator amenable, $(\ast)$ even splits: We obtain a (completely bounded) projection $P \!: H_c \to K_c$ which is also a left $\A$-module 
homomorphism. (The required splitting result can easily be proven by a more or less verbatim copy of the proof of its classical counterpart \cite[Theorem 2.3.13]{Run}.)
\par
The remainder of the proof is like in the classical case.
\par
For $f \in {\cal C}(\Omega)$, let $M_f$ denote the corresponding multiplication operator on $H$. The fact that $P$ is an $\A$-module homomorphism means that
\[
  M_f P = P M_f \qquad (f \in \A).
\] 
Since each $M_f$ is a normal operator with adjoint $M_{\bar{f}}$, the Fuglede--Putnam theorem implies that
\[
  M_{\bar{f}} P = P M_{\bar{f}} \qquad (f \in \A)
\] 
as well, and from the Stone--Weierstra{\ss} theorem, we conclude that
\[
  M_f P = P M_f \qquad (f \in {\cal C}(\Omega)).
\] 
Since $\A$ is unital, this implies that $K = H$. 
\par
Let $f \in {\cal C}(\Omega)$ be arbitrary. Then there is a sequence $( f_n )_{n=1}^\infty$ in $\A$ such that $\| f - f_n \|_2 \to 0$. Hence, we have
\[
  \left| \int_\Omega f \, d\mu \right| = \lim_{n \to \infty} \left| \int_\Omega (f - f_n) \, d\mu \right| \leq
  \lim_{n \to \infty} \int_\Omega |f - f_n| \, d|\mu| \leq \lim_{n \to \infty} |\mu|(\Omega)^\frac{1}{2} \| f - f_n \|_2 \to 0,
\]
which is impossible because $\mu \neq 0$.
\end{proof}
\begin{corollary}
The following are equivalent for a uniform algebra $\A$:
\begin{items}
\item $\min \A$ is operator amenable.
\item $\A$ is operator amenable for any operator space structure on $\A$ turning $\A$ into a completely contractive Banach algebra.
\item $\A$ is amenable.
\item $\A$ is a commutative $\cstar$-algebra.
\end{items}
\end{corollary}
\begin{proof}
(i) $\Longrightarrow$ (iv) is the assertion of the theorem, and (iv) $\Longrightarrow$ (iii) is well known.
\par
(iii) $\Longrightarrow$ (ii) $\Longrightarrow$ (i): The amenability of $\A$ is equivalent to the operator amenability of $\max \A$. Let $\A$ be equipped with any operator space structure
turning it into a completely contractive Banach algebra. Since
\[
  \max \A \stackrel{\id}{\longrightarrow} \A \stackrel{\id}{\longrightarrow} \min \A,
\]
are surjective completely contractive algebra homomorphisms, it follows from basic hereditary properties of operator amenability (\cite[Proposition 2.2]{Rua2}) that the operator amenability of $\max \A$
entails that of $\A$ and, in turn, that of $\min \A$.
\end{proof}
\vfill
\begin{tabbing} 
{\it Address\/}: \= Department of Mathematical and Statistical Sciences \\
                 \> University of Alberta \\
                 \> Edmonton, Alberta \\
                 \> Canada T6G 2G1 \\[\medskipamount]
{\it E-mail\/}:  \> {\tt vrunde@ualberta.ca}\\[\medskipamount]
{\it URL\/}:     \> {\tt http://www.math.ualberta.ca/$^\sim$runde/}
\end{tabbing}

\end{document}